\documentclass[12pt, twoside, leqno]{article}
 \pdfoutput=1
\usepackage{amsmath,amsthm}
\usepackage{amssymb,latexsym}
\usepackage{enumerate}
\usepackage{lscape}

\theoremstyle{definition}

\newtheorem*{xrem}{Remark}

\begin{document}

\baselineskip=17pt

\title{A counterexample to 'Algebraic function fields with small class number'}

\author{Claudio Stirpe\\
}

\date{November 2013}

\maketitle

\renewcommand{\thefootnote}{}



\renewcommand{\thefootnote}{\arabic{footnote}}
\setcounter{footnote}{0}

\begin{abstract}
Using class field theory I give an example of a function field of genus $4$ with class number one over the finite field $\mathbb{F}_2$.
In a previous paper (see \cite{lemaqu}, Section 2) the authors gave a proof of the nonexistence of such a function field.
This counterexample shows that their proof is wrong.
\end{abstract}

\subsection*{Introduction}

In \cite{st2}, Section 5 a list of pointless curves over $\mathbb{F}_2$ is given. 
One of such curves of genus $4$ has class number one contradicting a result given in \cite{lemaqu}, Section 2.
This counterexample is obtained with class field theory.
For notation and basic definitions on ray class fields the reader can see \cite{au00} or \cite{st1}.

\subsection*{The counterexample}
Let $K$ be the rational function field $\mathbb{F}_2(x)$. Let $\mathfrak{m}$ be the place $(x^4+x+1)$ of $K$ 
and let $S$ be the place $(x^7+x^4+1)$. We denote by $K_S^\mathfrak{m}$ the ray class field of conductor
$\mathfrak{m}$ such that $S$ is split in $K_S^\mathfrak{m}/K$. This extension is abelian and the degree is 
$7\cdot (2^4-1)=105$ (see \cite{au00}, Example 1.5). Consider the subextension $F/K$ of degree $5$. The constant field 
of $F$ is $\mathbb{F}_2$ and the genus is $4$ by the Hurwitz genus formula. We check that the class number is one.

In \cite{maqu}, Section 3, it is proved that a function field over $\mathbb{F}_2$ of genus $4$ has class number one if
and only if there is only one place of degree $4$ and no place of smaller degree.
We have to check that the places of $K$ of degree smaller than $4$ (except for $\mathfrak{m}$)
are inert in $F/K$. As in \cite{st1}, Example 4.2 and Example 4.3 we compute the Frobenius automorphism $Frob(P)$ 
in $Gal(K_S^\mathfrak{m}/K)$ for all unramified places $P$ with $deg(P)\leq 4$.

When $P=(x+1)$ then we define $z=\frac{(x+1)^7}{x^7+x^4+1}\in \hat{K}^*_P\subseteq J$.
By the Local Artin Map $z$ corresponds to $Frob(P)^{v_P(z)}=Frob(P)^7\in D(P)\subseteq Gal(K_S^\mathfrak{m}/K)$.
Then $P$ is split in $F/K$ if and only if $Frob(P)$ belongs to the subgroup $Gal(K_S^\mathfrak{m}/F)$ of
$Gal(K_S^\mathfrak{m}/K)$ of order $21$ (see \cite{st1}, Proposition 2.2).
But $z^3\not \equiv 1$ mod $x^4+x+1$ so the class of $z^3$ in the class group is not in $C_S^\mathfrak{m}$, 
the kernel of the Artin map (see
\cite{st1}, Definiton 3.6) and so $Frob(P)$ has not order $21$ in $Gal(K_S^\mathfrak{m}/K)$. It follows that $P$ is inert in $F/K$. 

The other unramified places of $K$ of degree smaller than $4$ are checked to be inert in a similar way.
In particular one can check the other rational places by considering $z=\frac{x^7}{x^7+x^4+1}$ and
$z=\frac{1}{x^7+x^4+1}$. Moreover for the places of degree two or three we can do similar computations with
$z=\frac{(x^2+x+1)^7}{(x^7+x^4+1)^2}$ or $z=\frac{(x^3+x+1)^7}{(x^7+x^4+1)^3}$
and $z=\frac{(x^3+x^2+1)^7}{(x^7+x^4+1)^3}$, respectively.
Finally for unramified places of degree four we have to consider
$z=\frac{(x^4+x^3+1)^7}{(x^7+x^4+1)^4}$ and $z=\frac{(x^4+x^3+x^2+x+1)^7}{(x^7+x^4+1)^4}$.

\begin{xrem}The choice of $S=(x^7+x^4+1)$ is not unique, in fact similar computations show that $S=(x^7+x^3+1)$
gives an other ray class field extension $K_S^\mathfrak{m}/K$ such that the unique subextension of degree five $F/K$
has no split places of degree smaller than $4$.
\end{xrem}

\small
\begin{center}{Claudio Stirpe, 
 E-mail: clast@inwind.it.}
\end{center}

\begin{thebibliography}{HD}

\baselineskip=17pt

\bibitem{au00} R. Auer, \emph{Ray class fields of global function fields with many rational places}, Acta Arithmetica 95, 97-122 (2000).
\bibitem{lemaqu} J. Leitzel, M. Madan, C. Queen \emph{Algebraic function fields with small class number}, Journal of Number Theory 7, 11-27 (1975).
\bibitem{maqu} M. Madan, C. Queen \emph{Algebraic function fields of class number one}, Acta Arithmetica 20, 423-432 (1972).
\bibitem{st1} C. Stirpe, \emph{An upper bound for the genus of a curve without points of small degree}, Phd Thesis at Universit\`a di Roma 'Sapienza' (2011), http://padis.uniroma1.it/bitstream/10805/1371/1/tesi.pdf.
\bibitem{st2} C. Stirpe \emph{An upper bound for the minimum genus of a curve without points of small degree}, Acta Arithmetica 160, 115-128  (2013).
\end {thebibliography}
\end{document}